\documentclass[11pt]{article}
\usepackage{graphicx}
\usepackage{color}
\usepackage{amsmath,amssymb,theorem}
\usepackage{times}
\usepackage{array}
\usepackage[english]{babel}
\usepackage{mathrsfs,color,comment}
\usepackage{a4wide}

\definecolor{wineRed}{rgb}{0.7,0,0.3}

\newtheorem{proposition}{Proposition}[section]

\newtheorem{remark}{Remark}[section]

\title{Some comments on using fractional derivative operators in modeling non-local diffusion processes}

\author{T.Namba${}^1$, P. Rybka${}^2$, V.R. Voller${}^3$\\
${}^1$ Nippon Steel \& Sumitomo Metal Corporation, Research \& Development\\
20-1 Shintomi, Futtsu, Chiba Prefecture
293-8511, Japan\\ e-mail: {\tt  namba.rb3.tokinaga@jp.nssmc.com}\\
${}^2$ Institute of Applied Mathematics and Mechanics,
Warsaw University\\ ul. Banacha 2, 07-097 Warsaw, Poland\\
fax: +(48 22) 554 4300, e-mail: {\tt rybka@mimuw.edu.pl}\\
${}^3$ Civil, Environmental, and Geo-  Engineering, University of Minnesota\\
500 Pillsbury Drive SE, Minneapolis, MN 55455, USA,\\
e-mail: {\tt  volle001@umn.edu}
}

\begin{document}
\maketitle
\begin{center}

\end{center}

\abstract
We start with a general governing equation for diffusion transport, written in a conserved form, in which the phenomenological flux laws can be constructed in a number of alternative ways. We pay particular attention to flux laws that can account for non-locality through space fractional derivative operators. The available results on the well posedness of the governing equations using such flux laws are discussed. 
A discrete control volume numerical solution of the general conserved governing equation is developed and a general discrete treatment of boundary conditions, independent of the particular choice of flux law, is presented. We use numerical solutions of various test problems to compare the operation and predictive ability of two discrete fractional diffusion flux laws  based on the Caputo (C) and Riemann-Liouville (RL) derivatives respectively. When compared with the C flux-law  we note that the RL flux law includes an additional term, that, in a phenomenological sense, acts as an apparent advection transport. Through our test solutions we show that, when compared to the performance of the C flux-law,  this extra term can lead to RL-flux law predictions that may be physically and mathematically unsound. We conclude, by proposing a parsimonious definition for a fractional derivative based flux law that removes the ambiguities associated with the selection between non-local flux laws based  on the RL and C fractional derivatives.    



\section{Introduction}

A non-local transport phenomena can be defined as one where the transport flux at given point depends not only on values at that point but also on values at remote locations. For example, a non-local diffusion flux can be defined as a weighted sum of potential gradients evaluated at points throughout the domain. Due to their non-local nature many works (e.g., see \cite{Zhang, Rina, VollerAHT} have been directed at using fractional derivative operators to construct non-local transport models in physical settings. 
Recent work \cite{Baeumer} studied the numerical solution of the diffusion equation  
\begin{equation}\label{nonc}
\partial_t u(x,t)=\kappa L^{1+\alpha} u(x,t) + f, \  \ 0<\alpha \le 1, \ \ x \in \Omega, \ t>0,
\end{equation} 
where $\Omega$ is the open interval $(0,1)$,  the operator $L^{1+\alpha}$ is a fractional derivative of order $1< 1+\alpha \le 2$, $\kappa$ is an appropriately scaled diffusivity, and $f$ is a source/sink term.  By numerically predicting the fate of an initial non-negative pulse $u(x,0) \geq 0$ on $\Omega$, the authors of \cite{Baeumer} investigate three alternatives for this operator, a Riemann-Liouville (RL$^{1+\alpha}$) derivative, a Caputo (C$^{1+\alpha}$) derivative, and a so called  \textit{Patie-Simon} (PS$^{1+\alpha}$) derivative. A major and important contribution from the work of \cite{Baeumer} is the development of a  comprehensive and consistent treatment of boundary conditions for numerical solution of the fractional heat conduction equation. They demonstrate, under a variety of boundary and initial conditions, that numerical solutions using the RL$^{1+\alpha}$ and PS$^{1+\alpha}$ derivatives, while not-identical, are always conservative. This, however, is not the case for the C$^{1+\alpha}$ derivative, under some boundary and initial conditions this operator may not be conservative. The authors of \cite{Baeumer} demonstrate this by considering a problem with fixed boundary values, $u(0,t)=u(1,t)=0$, presenting a numerical solution, based on the C$^{1+\alpha}$ derivative, where an initial non-negative pulse $u(x,0)$ produces at a later time solutions $u(x,t)$ with negative values.        

An alternative form for the fractional diffusion equation is 
\begin{equation}\label{con}
\partial_t u(x,t)=\partial_x\left(\kappa L^{\alpha} u(x,t)\right) + f(x,t) , \  \ 0<\alpha \le 1, \ \ x \in \Omega,\ t>0,
\end{equation} 
where we interpret the right-hand side as the conserved statement the \em divergence of a flux \em $q= L^{\alpha} u(x,t)$. We note that in using the equation form suggested by \cite{Baeumer}, i.e., eq. (\ref{nonc}), conservation will only by ensured if $\partial_x\left(L^{\alpha}u(x,t)\right) \equiv L^{\alpha+1}u(x,t)$. This is true in the case 
the  operator $L^{\alpha}$ is the RL$^{\alpha}$ derivative of order $0 < \alpha \le 1$. For Caputo derivatives, however, the condition that $\partial_x\left(C^{\alpha}u(x,t)\right) \equiv C^{\alpha+1}u(x,t)$, only holds if $\partial_x u(0,t)=0$ (see Remark 6.5 in \cite{Baeumer} and Subsection 3.1. below).  This explains, why, in some cases use of the order ($1 < 1+\alpha \le 2)$ C$^{\alpha+1}$  fractional derivative in eq. (\ref{nonc}) is not  conservative. Indeed, we note that a core advantage of the  PS$^{\alpha+1}$ derivative, adopted by \cite{Baeumer}, is that it by its construction,  $\partial_x\left(\text{C}^{\alpha}u(x,t)\right) \equiv \text{PS}^{\alpha+1}u(x,t)$ and thus, when used in eq.(\ref{nonc}) conservation is enforced.

Our task here is to study the space fractional diffusion equation in $\Omega$, under a variety of  flux definitions. In a departure from \cite{Baeumer}, we will base our analysis and solutions on the conserved form of the diffusion  transport equation, eq. (\ref{con}). Our investigation generates a number of insights and results related to the solution of fractional diffusion equations in bounded domains and also highlights possible advantages in using the Caputo derivative in seeking discrete (numerical) solutions.


We start with a general conserved statement for a transient diffusion
transport written in terms of a general flux $q$. On specific
definition of the flux we arrive at three alternative diffusion
transport models, the classic Fourier second law ($q\sim \partial_x
u$), a Caputo model ($q \sim \mbox{C}^{\alpha}$), and a
Riemann-Liouville model ($q\sim \mbox{RL}^{\alpha}$).

To lay a foundation for development of numerical solutions for the
general conserved model  we present some analytical results related to
the well posedness of fractional flux models.

We gather here in Section 3 the currently available results, see \cite{BJDE},
\cite{Barxive}, \cite{NaRy}, \cite{Rys}, on the well-posedness of (\ref{nonc}) or 
(\ref{con}) augmented with an array of boundary conditions and initial data.
It turns out that we can consider the data from the following list: Dirichlet 
conditions, fixed flux, fixed slope conditions and mixed ones.

Another important issue is the validity of the maximum principle. It holds for some problems written in a
conservative form. 
We present not only  a rigorous argument, but also discuss  the consequence of this principle  on the invariance of eq. (\ref{con}) under a re-scaling; i.e., regardless of the measurement scale used the actual value of the potential $u$ (e.g., temperature) should always be the same. We show, by example, that eq. (\ref{con}) with the flux given in terms of the Riemann-Liouville derivative fails to satisfy this invariance. On the the other hand this scale invariance and the maximum principle hold if we use the Caputo derivative in eq. (\ref{con}). We offer the explanation in Section
\ref{six} and further comments in Section \ref{seven}.

A thorough discussion of the boundary value problem (\ref{con}) is not complete without addressing the issue of the boundary conditions. This is done in Section \ref{three}, where we ask the question whether or not we can evaluate $u$, the flux or the derivative of $u$ at the boundary of $\Omega$ for $u$ a solution to (\ref{con}). In particular, we note that the condition $\partial^\alpha_{x_C}u(0) = 0$ is automatically satisfied for sufficiently regular $u$, see Proposition \ref{flux}. We explain it in  detail in Section 3.1.
This  section is a good place to  to discuss the 
correctness of the  
definition of $\partial^\alpha_{x_C}$ and $\partial^\alpha_{x_{RL}}$.

Following this theory, we develop an explicit control volume discretization for the general conserved model, identifying a range of possible boundary condition treatments, the operation of which are essentially invariant to the precise definition of the flux used. Then, taking a similar approach to the one adopted in \cite{Baeumer} we study the fate of an initial non-negative pulse in $\Omega$ with both the Caputo and Riemann-Liouville models. As expected these test gives identical results to those provided in \cite{Baeumer} and illustrate the predictive differences between using Caputo and Riemann-Liouville flux models. However, due to the fact that our analysis is based on a conserved treatment we are able to write down explicit discrete expressions for the Riemann-Liouville and Caputo fluxes. These expressions show that difference between the Caputo and Riemann-Liouville flux models reduce to the addition of a correction term to the Riemann-Liouville flux with the form of an apparent advection transport. This observation provides us with a direct and clear explanation of the predictive differences between Caputo and Riemann-Liouville fluxes. In addition it also allows us to identify situations where the Riemann-Liouville flux definition may lead to predictions which are physically and mathematically unsound. To conclude we suggest an alternative fractional flux law that essentially negates the differences between the Caputo and Riemann-Liouville flux models.
 
\section{The conserved diffusion transport equation}

The transient diffusion equation, on the interval $\Omega = (0,1)$, written in the conserved form is
\begin{equation}\label{eqbal}
\partial_t u =-\partial_x q+f
\end{equation}
where $q$ is the diffusion flux. 

We consider three separate definitions for this flux:

\noindent $\bullet$
\textbf{The Fourier Law}
Assuming a unit diffusivity $\kappa =1 $, the classic flux definition is the Fourier law \cite{Joe}
\begin{equation}\label{F}
q(x,t)=-\partial_x u(x,t)	\equiv -u_x(x,t)
\end{equation}
stating that the flux is proportional to the potential gradient.

\noindent $\bullet$ \textbf{Riemann-Liouville flux}
The Fourier law is \emph{local}, i.e., calculating the flux at a specified point only requires knowledge of the value of the gradient at that point. In systems where we have a range of heterogeneity length scales \cite{Metzler}, however, the evaluation of the flux at a specified point may also be controlled by conditions remote from the point of interest. In such cases a candidate model for the flux can be formed from a fractional derivative \cite{Pod} of the potential which, as noted in our introduction, can be regarded as a weighted sum of potential gradients evaluated at points throughout the domain. One choice for this flux model is the Riemann-Liouville (RL) fractional derivative of order $0<\alpha \le 1$, providing the following flux definition. 
\begin{equation*}
q(x,t)=-\partial_{x_{RL}}^{\alpha} u(x,t)
\end{equation*}
where the derivative $\partial_{x_{RL}}^{\alpha} u$ is defined below,
\begin{equation}\label{RL}
\partial_{x_{RL}}^{\alpha}u(x,t)=\frac{1}{\Gamma(1-\alpha)} \frac{\partial}{\partial x} \int_0^x(x-\xi)^{-\alpha} \  u(\xi,t) \, d \xi	
\end{equation}
\noindent $\bullet$  \textbf{Caputo flux}
As an alternative to the Riemann-Liouville (RL) fractional  flux we
can use a  Caputo (C) fractional operator of order $0<\alpha \le 1$,
i.e.
$$
q(x,t)=-\partial_{x_{C}}^{\alpha} u(x,t),
$$
where the derivative $\partial_{x_{C}}^{\alpha} u$ is given by the
following formula,
\begin{equation}\label{C}
\partial_{x_{C}}^{\alpha} u(x,t) =\frac{1}{\Gamma(1-\alpha)}\int_0^x(x-\xi)^{-\alpha}\ 
\frac{\partial u}{\partial \xi}(\xi,t) \, d \xi.	
\end{equation}
We recall the important properties of the C derivative \cite{Pod, Mark},  which we will use:\\

\bigskip\noindent
\textbf{Fact 1.}  The  C and RL  derivatives of order $0<\alpha \le 1$ of the function
  $u^*(x,t)$  are equivalent if the left-hand boundary data of the
  function vanishes i.e., if $u^*(0,t)=0$. However, we must make sure
  that the  $u^*(0,t)$ is well-defined, i.e. $u^*$ has a trace at the
  boundary, see the discussion in \S 3.1. below.
\begin{proposition}
 The order $0<\alpha \le 1$ C derivative of a constant $a$ is zero, i.e., $\partial_{x_{C}}^{\alpha} a =0$.
\end{proposition}
This fact is obvious from the definition of the integer derivative, however, it does not hold for the Riemann-Liouville derivative.

\bigskip\noindent{\textbf{Fact 2.}  From the definition of the RL derivative, eq. (\ref{RL}), it immediately   follows that the 1st order derivative of the order $0<\alpha \le 1$ C derivative is the $0<\alpha \le 1$ RL derivative of the Fourier flux, i.e., 
\begin{equation}\label{r7}
\partial_x \left(\partial^\alpha_{x_C} u(x,t)\right)=\partial^\alpha_{x_{RL}} \left(\partial_{x} u(x,t)\right).
\end{equation}
This property allow us to exploit desirable properties of the RL derivative in analyzing the well posedness of the conserved equation with a C flux definition. 
 

\section{Well posedness of diffusion equation solutions based on a fractional flux models}\label{three}

\subsection{On the definitions of fractional derivatives}
The first issue to be clarified before tackling the differential equations is the well-posedness of the operators $\partial^{\alpha}_{x_{RL}}$ and $\partial^{\alpha}_{x_C}$ for $\alpha\in (0,1)$. It is convenient to use for this purpose the notion of the weak derivative and Sobolev spaces. The space $W^{k,p}(\Omega)$, $k\in \mathbb{N}$, consists of functions whose weak derivatives up to order $k$ are in $L^p(\Omega)$,  $p\ge 1$ and $H^s(\Omega)$, $s>0$ is the fractional Sobolev space, we note that for $k\in \mathbb{N}$ we have $W^{k,2}(\Omega) = H^k(\Omega)$. A good source of information about the Sobolev spaces is \cite{adams}.

If we look at the definition (\ref{C}) of the Caputo derivative, then we see that $u$ must have a weak derivative, which is at least integrable, i.e. $u \in W^{1,1}(\Omega)$. 
Let us apply the point of view of Functional Analysis on 
$\partial^{\alpha}_{x_{RL}}$ and consider this as an unbounded operator on $L^2(\Omega)$, i.e. $\partial^{\alpha}_{x_{RL}} : 
D( \partial^{\alpha}_{x_{RL}}) \subset L^2(\Omega) \to L^2(\Omega)$. It turns out that this problem was studied in \cite{Yamamoto}. Roughly speaking his result may be conveyed as follows, see \cite[Theorem 2.1]{Yamamoto} for details, 
$$
D( \partial^{\alpha}_{x_{RL}}) \subset H^{\alpha}(\Omega).
$$
It is now clear that the definition of
$\partial^{\alpha}_{x_C}$ requires more stringent conditions on
$u$, than $\partial^{\alpha}_{x_{RL}}$. This is also clear from the equivalent definition of
$\partial^{\alpha}_{x_C}$, mentioned in Fact 1, i.e.
$$
\partial^{\alpha}_{x_C} u  =  \partial^{\alpha}_{x_{RL}}(u - u(0)).
$$
On the one hand we relax the requirements on $u$, i.e. $u \in H^\alpha$ will suffice, but the requirement that the trace of $u(0)$ is well-defined restricts the range of $\alpha$ to $\alpha >1/2$.

Even though, the operator $\partial^\alpha_{x_C}$ may be properly defined for $\alpha\le \frac12$, we will restrict our attention to the range $(\frac12, 1]$ for the sake of technical simplicity.

Let us also comment on the zero-flux boundary condition, i.e., $q=-\partial^{\alpha}_{x_C}=0$.
In \cite{Rys} the
author studied the domain in 
$ L^2(\Omega)$
of $\frac\partial{\partial x} \partial^{\alpha}_{x_C}$ with the following boundary conditions,
\begin{equation}\label{kropa}
u_x(0,t) = 0 = u(l,t),\qquad t>0.
\end{equation}
It is shown in  \cite{Rys} that 
\begin{equation}\label{rd}
 D(\frac\partial{\partial x} \partial^{\alpha}_{x_C})\subset H^{1+\alpha}.
\end{equation}
Keeping this in mind we note that frequently 
the zero-flux condition is automatically satisfied at $x=0$. Indeed, we have:

\begin{proposition}\label{flux}
If $u\in H^{1+\alpha}(\Omega)$ and $\alpha>1/2$, then
\begin{equation}\label{rf}
 \partial^{\alpha}_{x_C} u(0) = 0.
\end{equation}
\end{proposition}
{\it Proof.}
Indeed, if $\alpha > 1/2$, then for $u \in H^{1+\alpha}$ we have that $u_x  \in H^\alpha$, i.e. $u_x$ is continuous in $\bar\Omega$, hence $u_x$ is bounded by a constant $M>0$. Thus,
$$
|\partial^{\alpha}_{x_C} u (x) | =  c  \left| \int_0^x \frac{u_x(s)}{(x-s)^\alpha} \,ds \right|
\le c \frac M{1-\alpha} x^{1-\alpha} \to 0\qquad\hbox{as } x\to 0.\eqno\Box
$$

\begin{remark}
 We can  draw a similar conclusion by using the H\"older inequality if we happen to know that $u \in W^{1, 1/(1-\alpha)}$
for $\alpha\in (0,1)$.
\end{remark}

\subsection{On the well-posedness of a diffusion equation based on a fractional flux models}
The literature on the rigorous analytical treatment of
\begin{equation}\label{nr0}
 u_t = \partial^{1+\alpha}_{x} u +f\qquad\hbox{in } Q_T:= (0,1)\times(0,T)
\end{equation}
or 
\begin{equation}\label{nr1}
 u_t = \frac \partial {\partial x} (\partial^\alpha_{x_C} u) +f\qquad\hbox{in } Q_T:= (0,1)\times(0,T)
\end{equation}
for $\alpha\in (0,1)$, with the initial condition
\begin{equation}\label{nr3}
 u(x,0) = u_0(x) \equiv g(x)
\end{equation}
and appropriate boundary conditions is rather scant. Here, $\partial^{1+\alpha}_{x} u$
denotes the Riemann-Liouville or Caputo derivative. 
These problems we studied in \cite{BJDE}, \cite{Barxive}, \cite{NaRy}, \cite{Rys} and \cite{Rys2}.

On the one hand, the authors of \cite{BJDE} address both problems for a variety of boundary conditions, i.e. they study:\\
(1) the absorbing, i.e. Dirichlet, boundary conditions; \\
(2) the zero flux, i.e. reflective boundary conditions,
$$
\partial^{\alpha}_{x} u (x,t)= 0
$$
at $x \in \partial \Omega$;\\
(3) the zero Neumann, i.e. zero slope boundary conditions;\\
(4) mixed boundary conditions, see (\ref{kropa}). However, they are not dealt with explicitely.\\
On the other hand the authors of   \cite{Barxive} concentrate only on the absorbing boundary conditions.


The papers  \cite{BJDE} and  \cite{Barxive} are written from the probabilistic point of view, treating the right hand sides of equations (\ref{nr0}) and (\ref{nr1}) as generators of strongly continuous semigroups. The papers \cite{NaRy} and \cite{Rys} offer the perspective of the Partial Differential Equations on (\ref{nr1}). In \cite{NaRy} we deal with the  general 
Dirchlet boundary data, however here, for the sake of simplicity we consider
\begin{equation}\label{nr2}
 u(0,t) = 0 = u(1,t)\qquad\hbox{for } t\in(0,T).
\end{equation}
In \cite{Rys, Rys2} equations (\ref{nr1}), (\ref{nr3}) are considered with  mixed boundary conditions,
\begin{equation}\label{nr2m}
 u_x(0,t) = 0 = u(1,t)\qquad\hbox{for } t\in(0,T).
\end{equation}
We argued in Proposition \ref{flux} that the  flux condition is
automatically satisfied at $x=0$, if $u\in H^{1+\alpha}$ and $\alpha>1/2$.

Interestingly, 
the methods used in both papers, \cite{NaRy, Rys} are not capable to deal with the case of zero slope at the boundary, 
$$
u_x(0,t) = 0 = u_x(1,t)\qquad\hbox{for } t\in(0,T).
$$
Both papers, \cite{NaRy} and \cite{Rys}, offer well-posedness
results. We mention that problem (\ref{nr1}), (\ref{nr3}), (\ref{nr2}) is studied in \cite{BJDE} from a completely different perspective of the $C_0$-semigroups. At the same time problem (\ref{nr1}), (\ref{nr3}), (\ref{nr2m}) is not studied there.

We will now briefly comment on the
content of \cite{NaRy} and \cite{Rys, Rys2}.

The purpose of \cite{NaRy} is to introduce the notion of viscosity solutions, which is based on the maximum principle. It is well-known that such an approach is quite successful for evolution equations involving time fractional derivatives or non-local spacial operators, see \cite{Barles, GigaNamba, Caffarelli}. 

We claim that our notion of viscosity solution is correct in the sense that one can show existence and  uniqueness of solutions to the problem (\ref{nr1}-\ref{nr3}), as well as the continuous dependence of solutions upon the data. Indeed, in \cite[Theorem 5.6]{NaRy} we prove existence and uniqueness of solutions to  (\ref{nr1}-\ref{nr3}) provided that functions $f$ and $g$ are continuous and bounded. To achieve this goal we apply the classical Perron method. 
Moreover, the continuous dependence of solutions upon the data $u_0$ and $f$ are shown there too, see \cite[Proposition 6.2]{NaRy}. 
These results are based on  comparison principle, see \cite[Theorem 4.1]{NaRy}. Its simple statement may go like this.

\begin{proposition}
 Let us suppose that $u$ is a (viscosity) solution to (\ref{nr1}), (\ref{nr3}), (\ref{nr2}), where $f\equiv 0$, $g\in C([0,1])$ and $g(0) = 0= g(1)$. Finally, we assume that $g(x) \ge 0$ for all $x\in[0,1]$. Then, 
 $$
 \max_{\{x\in \bar\Omega} g(x) \ge u(x,t)\ge 0\quad\hbox{for all }(x,t) \in Q.
 $$
\end{proposition}
In particular, this Proposition implies uniqueness of viscosity solutions to  (\ref{nr1}-\ref{nr2}).

A completely different approach to the question of existence of solutions to (\ref{nr1}), (\ref{nr3}), (\ref{nr2m}) is applied in \cite{Rys}. It is based on the theory of semigroups. Actually, Ryszewska shows in \cite[Theorem 2]{Rys} the operator $\frac{\partial}{\partial x}\partial^\alpha_{x_C}$ defined on an appropriately chosen domain in $L^2(0,1)$ generates an analytic semigroup. Thus, we immediately obtain existence, uniqueness and regularity of solutions, see \cite[Theorem 3]{Rys}. Moreover, we also have the continuous dependence of solutions upon the data.

Interestingly, a comparison principle holds too, see \cite[Lemma 9]{Rys2}.

\begin{proposition}
Let us suppose that $u$ is a (semigroup) solution to (\ref{nr1}), (\ref{nr3}), (\ref{nr2m}), where $f\equiv 0$, $g\in C([0,1])$ and $g_x(0) = 0= g(1)$. Finally, we assume that $g(x) \ge 0$ for all $x\in[0,1]$. Then, 
 $$
 u(x,t)\ge 0\quad\hbox{for all }(x,t) \in Q.
 $$
\end{proposition}


We emphasize that in \cite{NaRy} and \cite{Rys}
the maximum principle is stated  explicitly in the above Propositions. Of course, the maximum principle implies the minimum principle, after the sign reversal.

The authors of \cite{BJDE} and  \cite{Barxive} are aware that the maximum principle may be violated for  (\ref{nr0}), when they write that positivity of solutions is not preserved. 
In particular, the authors  of \cite{Baeumer}  found a numerical example of initial condition for
\begin{equation}\label{nr4}
u_t = \partial_{x_C}^{1+\alpha} u, \qquad(x,t)\in Q,
\end{equation}
where $\alpha\in(0,1)$
with boundary conditions (\ref{nr2}) and initial data (\ref{nr3}),
such that the  solution was negative at later times.

\section{A discrete form of the conserved equation}

We will approximate the conserved transport equation eq. (\ref{con}) using a control volume scheme \cite{Pat,VollerBook}. The first step is to lay down a grid of $n+1$ equally spaced ($\Delta x =1/n$) node points, see Fig \ref{fig1}. Around the internal nodes $i=1,2,..,n-1$, we construct control volumes of length $\Delta x$, for a given node point located at $i\Delta x$,  the $in$ and $out$ faces of the control volume located at positions $x=(i-0.5)\Delta x$ and $x=(i+0.5)\Delta x$ respectively. Half control volumes (length $\Delta x/2$) are associated with the end-point nodes $i=0$ and $i=n$. The control volume faces associated with the node $i=0$ are located at $x=0$ and $x=0.5\Delta x$ and those associated with node $n$ are located at $x=1-\Delta x/2$ and $x=1$. At any internal node a conserved discrete version of eq. (\ref{eqbal}) can be formed by taking backward difference in time and central difference in space approximations to arrive at, neglecting source and sink terms,
\begin{figure}[h]
\centering
\includegraphics[width=10cm]{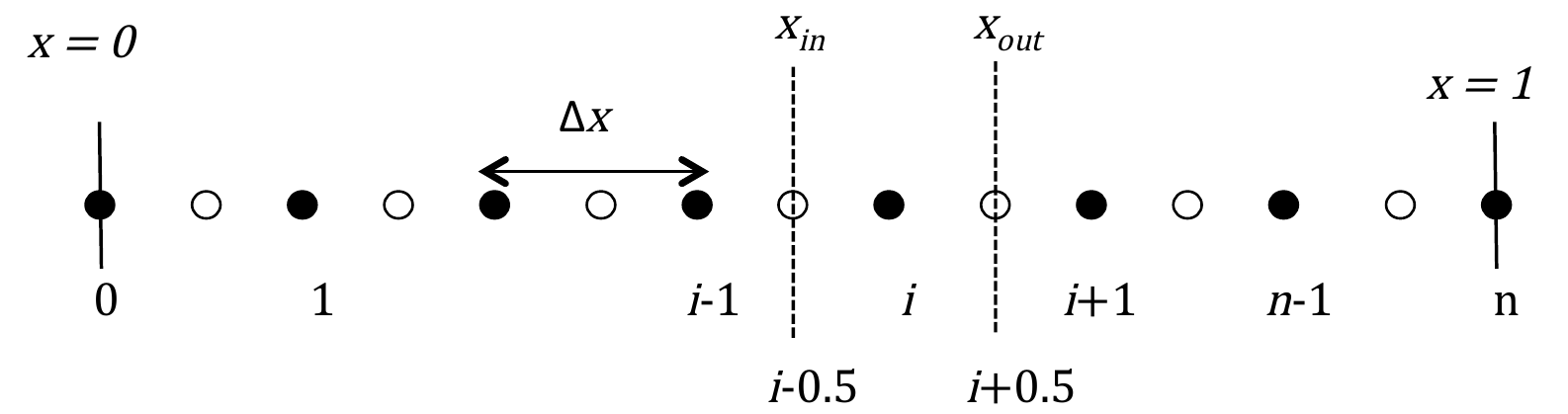}
\caption{Control volume grid}
\label{fig1}
\end{figure}
\begin{equation}\label{eqdis}
u_i^{k+1}=u_i^k+\frac{\Delta t}{\Delta x}\left(q_{i-0.5}^k-q_{i+0.5}^k \right )
\end{equation}
where the super-script $k$ indicates evaluation at time $t=k\Delta t$ and the sub-script $i$ indicates evaluation at $x=i\Delta x$. Note that the fluxes are evaluated at the $in$ and $out$   control volume faces.

Here we will consider two possible boundary conditions for the discrete balance eq. (\ref{eqdis}). 

%
\paragraph{Fixed Value}

{When we have a fixed value ( i.e., a Dirichlet condition) of the potential $u$ at the end points ($u_{left}$ or $u_{right}$)  of the domain we simply set  
\begin{equation}\label{left}
	u_0^k=u_{left}
\end{equation}                
if the condition is at $x=0$, or
\begin{equation}\label{right}
	u_0^k=u_{right}
\end{equation}                
if the condition is at $x = 1$. }

We note that when a fixed boundary value  of 0  is applied, the associated boundary is said to be absorbing. 

\paragraph{Fixed Flux}

The second boundary condition is a fixed flux condition $q|_{x=0} = q_{left}$ or $q|_{x = 1} = q_{right}$. If applied at $x=0$, such a condition will result in the following discrete scheme in the 1/2 control volume associates with node 0,
\begin{equation}\label{fixleft}
	u_0^{k+1}=u_0^k+2\frac{\Delta t}{\Delta x}\left(q_{left}-q_{0.5}^k \right ),
\end{equation}
where the factor 2 is due to the 1/2 volume size and  the subscript $0.5$ implies evaluation at the face between nodes 0 and 1. In a similar manner the fixed flux condition is applied at $x=1$ we modify the scheme at node $n$
\begin{equation}\label{fixright}
u_n^{k+1}=u_n^k+2\frac{\Delta t}{\Delta x}\left(q_{n-0.5}^k-q_{right} \right),
\end{equation}
where the subscript $n-0.5$ indicates evaluation at the face between nodes $n-1$ and $n$.  

A fixed flux case of particular interest is when the specified flux is
zero, i.e, $q_{left}=0$ and/or $q_{right}=0$, this condition is
referred to in \cite{Baeumer} as a reflective boundary condition. A
closely related condition is to use a Neumann condition and set the
slopes at the boundaries to zero, i.e., $u_x(0,t)=0$ and/or
$u_x(1,t)=0$. Depending on the choice of the flux model and the
location of the boundary, zero-slope  conditions may or may not be
equivalent to zero flux conditions, see also \S 3.1.      

%

\section{Specific discrete flux definitions}
The convenience of eq. (\ref{eqdis}), and its boundary treatments eqs (\ref{left}) to (\ref{fixright}), is that, due to the conservative form, its operation will be invariant under any appropriate flux definition. Thus making the numerical solution simply requires the appropriate construction of the discrete fluxes at the control volume faces. 

\subsection{Fourier flux}
When we have the Fourier flux defined in eq. (\ref{F}) the fluxes and the $in$ and $out$ faces of the i-th control volume can be discretely approximated by the central finite difference scheme,
\begin{equation}\label{fflux}
\begin{aligned}
q^k_{i-0.5_{F}}&\approx\frac{u^k_{i-1}-u^k_{i}}{\Delta x}, \\
\\
q^k_{i+0.5_{F}}&\approx\frac{u^k_{i}-u^k_{i+1}}{\Delta x}.	
\end{aligned}
\end{equation}
%
%

\subsection{Riemann-Liouville fractional flux}

At an internal control volume face, we can evaluate a discrete form of the Riemann-Liouville flux, defined in eq. (\ref{RL}) flux, with the  $\tfrac{1}{2}$ shifted Gr\"unwald weights \cite{Mark,Dan}, e.g., at the $in$ and $out$ face of the volume $i$ (see Fig. \ref{fig1}) 
\begin{equation}\label{eqgflux}
\begin{aligned}
q^k_{i-0.5_{RL}}&=-\left(\partial_{x_{RL}}^{\alpha} u\right)_{x=(i-0.5)\Delta x} \approx-\Delta x^{-\alpha}\sum_{j=0}^{i} g_j u^k_{i-j}, \\
\\
q^k_{i+0.5_{RL}}&=-\left(\partial_{x_{RL}}^{\alpha} u\right)_{x=(i+0.5)\Delta x} \approx-\Delta x^{-\alpha}\sum_{j=0}^{i+1} g_j u^k_{i+1-j}.
\end{aligned}
\end{equation}
%
%
%
where the Gr\"unwald weights are given by
\begin{equation}\label{gweights}
\begin{aligned}
g_0 &=1 , \\ 
g_i &=\frac{i-1-\alpha}{i}g_{i-1},\quad i =1,2, \ldots .
\end{aligned}
\end{equation}
Through a combination of addition and subtractions \cite{VollerAHT, Dan} we can rewrite the face flux evaluations in eq. (\ref{eqgflux}) in the alternative form,
\begin{equation}\label{eqrlwflux}
\begin{aligned}
q^k_{i-0.5_{RL}}&\approx\sum_{j=0}^{i-1} W_j\frac{u^k_{i-1-j}-u^k_{i-j}}{\Delta x}-\frac{W_{i}}{\Delta x} u^k_0, \\
\\
q^k_{i+0.5_{RL}}&\approx\sum_{j=0}^{i} W_j\frac{u^k_{i-j}-u^k_{i+1-j}}{\Delta x}-\frac{W_{i+1}}{\Delta x} u^k_0,
\end{aligned}
\end{equation}
where the weights $W_j \ge 0$, decreasing in value as $j$ increases, are given by 
\begin{equation}\label{weights}
W_j=\Delta x^{1-\alpha}\sum_{m=0}^j g_m .
\end{equation} 
Thus, apart from the single initial data term, we see, with  reference to eq. (\ref{eqgflux}), that the fractional flux, at a particular internal volume face, can be expressed as a weighted sum of the central difference approximations of the Fourier fluxes across the internal control volume faces, at and to the left of the face of interest.  Also note that  as $\alpha \rightarrow 1$, $W_0\rightarrow 1$, $W_{j>0} \rightarrow 0$ and,  as required, we recover the Fourier flux definitions from eq. (\ref{eqrlwflux}). Further, we note, that when these flux calculations are directly substituted into eq. (\ref{eqdis}) and its associated boundary treatments (eqs (\ref{left})--(\ref{fixright})), we arrive at schemes that are essentially identical to those presented in \cite{Baeumer} for the RL$^{\alpha+1}$ operator; the only difference is that, in the fixed flux case \cite{Baeumer} neglects the factor of 2 required to account for the half control volumes at the discrete domain boundaries. While the algebra may show that the resulting scheme are the same, we feel that the physical interpretation and general boundary treatment derived from using the conserved heat balance provides worthwhile benefit.             

\subsection{Caputo fractional flux}

Assuming that the boundary data on the left is always well defined, we can construct  the variable $u^*(x,t)=u(x,t)-u(0,t)$, which takes a zero value on the left boundary, i.e., $u^*(0,t)=0$. In this way, on noting that 
\begin{equation*}
u^{*k}_{i-1-j}-u^{*k}_{i-j}=(u^{k}_{i-1-j}-u^{k}_0)-(u^{k}_{i-j}-u^{k}_{0})=u^{k}_{i-1-j}-u^{k}_{i-j},
\end{equation*}
it follow from (\ref{eqrlwflux}), that the RL derivative of $u^*$, at location $i-0.5$ say, can be approximated as
    
\begin{equation}\label{derofustar}
\left(\partial_{x_{RL}}^{\alpha} u^*\right)_{x=(i-0.5)\Delta x}\approx-\sum_{j=0}^{i-1} W_j\frac{u^{k}_{i-1-j}-u^{k}_{i-j}}{\Delta x}. 
\end{equation}
We also know, from the properties of the C derivative listed below eq.(\ref{C}), that with zero left boundary data the C and RL derivative coincide, i.e., $\partial^{\alpha}_{x_{C}} u^*=\partial^{\alpha}_{x_{RL}} u^*$, and that, since at a given time $t$, functions $u^*(\cdot,t)$ and $u(\cdot,t)$ differ by a constant, then
$\partial^{\alpha}_{x_{C}} u^*=\partial^{\alpha}_{x_{C}} u$. From this it follows that $\partial^{\alpha}_{x_{C}} u= \partial^{\alpha}_{x_{RL}} u^*$, thus, from eq.(\ref{derofustar}), we can approximate the C-fluxes at the control volume faces as 
\begin{equation}\label{eqcwflux}
\begin{aligned}
q^k_{i-0.5_{C}}&\approx\sum_{j=0}^{i-1} W_j\frac{u^k_{i-1-j}-u^k_{i-j}}{\Delta x}, \\
\\
q^k_{i+0.5_{C}}&\approx\sum_{j=0}^{i} W_j\frac{u^k_{i-j}-u^k_{i+1-j}}{\Delta x}
\end{aligned}
\end{equation}
i.e., the C flux at a given control volume face is simply the weighted sum of the Fourier faces fluxes at and to the left.   


%
%
 %
%

Again we note that the scheme resulting on using the flux definitions of eq. (\ref{eqcwflux}) in (\ref{eqdis}) and its associated boundary treatments (eqs (\ref{left})--(\ref{fixright}))  is algebraically equivalent (apart from the boundary factor of 2) to the scheme derived by \cite{Baeumer} on using the  PS$^{\alpha+1}$ operator. In this case, however, we recognize  that the direct and uncorrected connection between the fractional flux and the weighted sum of Fourier integer fluxes is a \em property that should provide a significant level of robustness in using C derivatives in a non-local model of the flux.\em 
 
\section{Physical consequences of the difference between the discrete RL and C flux models}\label{six}

\subsection{Apparent advection}
We note that the difference between the discrete RL (eq. (\ref{eqrlwflux})) and C (eq. (\ref{eqcwflux})) fluxes is the extra term $(W_i/\Delta x)u_0^k$, involving the left hand boundary value, $u(0,t)\equiv u_0^k$, that appears in the RL flux. Here our intention is to provide a physical description for this term. 

We start by noting that, since it is purely comprised of a weighted sum of potential gradients,  in a phenomenological sense, the C flux  acts in the same way as the classic Fourier flux, eq. (\ref{eqgflux}). By this we mean that, if an initial pulse is provided within $(0,1)$, and  no-flux conditions are enforced at both boundaries, the C flux will continue to transport the conserved quantity, to both the left or the right, until all the local slopes, at all points in the domain, are zero. In all cases, in the long time limit, any initial profile for $u(x,0)$, will be averaged 
out to a  horizontal line, $u(x,t)=a$ . With this in mind, we will refer to the weighted sum, the sole term in the C flux in eq. (\ref{eqcwflux}), as the ``non-local diffusion". In contrast, faced  with the same no-flux boundary conditions on $(0,1)$, the  RL flux in  eq. (\ref{eqrlwflux}) will only return a horizontal line (constant $u$ value)  in the special   case that the initial data is $u(x,t)=0$. In the same manner as the C flux,  the ``non-local diffusion'', the first term in RL flux, will, depending on the sign of the weighted slope  at a given point, move the conserved quantity to the left or right in an  attempt to reduce ``peaks'' and fill in ``valleys''. This effort, however, is offset by the extra term in the RL flux. Noting that fact that the weights $W$, constructed for the flux definition, are positive this extra term will always transport the conserved quantity in one direction, towards the left boundary if the boundary data $u(0,t)>0$ and towards the right if $u(0,t)<0$.  In this way, from a phenomenological view point, we can consider this term in a similar manner to an  advection transport term. Thus we will refer to the extra term in the discrete RL flux (the second term in eq. (\ref{eqrlwflux})) as an ``apparent advection''. In essence, we are claiming that the operation of the discrete C flux definition eq. (\ref{eqcwflux}) mimics the  behavior of a diffusion transport process and that the discrete RL flux definition eq. (\ref{eqrlwflux}) mimics the behavior of a diffusive-advection process. 

In the following example problems, tracking, under given boundary conditions,  the fate of an initial $u(x,0)$ profile in (0,1), fully illustrates how the "apparent advection" term in the RL flux can prove problematic  in matching its predictions to physically observed behaviors. 

\subsection{Effect of apparent convection with reflective boundary conditions}
First we consider the fate of the initial pulse in $(0,1)$,
\begin{equation}
  u(x,0) =
    \begin{cases}
      0 & \text{if $x<0.3$},\\
       25x-7.5 & \text{if $x \ge 0.3$ and $x<0.5$},\\
      -25x+17.5 & \text{if $x \ge 0.5$ and $x<0.7$},\\
      0 &\text{if $x \ge 0.7$}\\
    \end{cases}       
\end{equation}
with zero flux (reflective) boundary conditions at $x=0$ and $x=1$; conditions that will ensure that the quantity in the initial pulse will be maintained through time.  Results of numerical calculations of this problem ($\alpha=0.5$) using, in turn, the flux approximations in eq.  (\ref{eqrlwflux}) (RL)  and eq. (\ref{eqcwflux}) (C), along with a time step of $\Delta t =0.0005$ and space step of $\Delta x=1/n, n=100$, are shown in Fig. \ref{fig2}.
\begin{figure}[h]
\centering
\includegraphics[width=12cm]{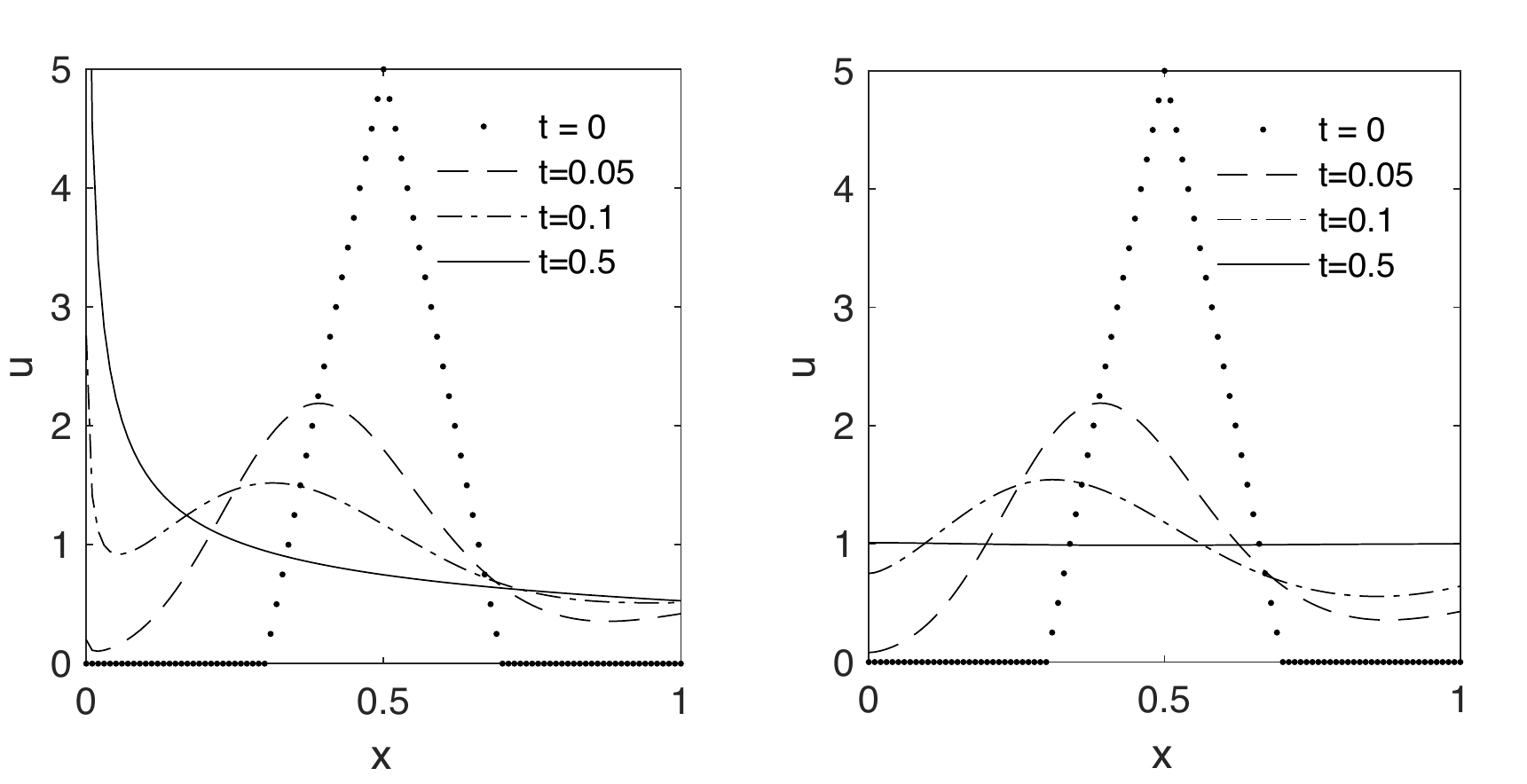}
\caption{Fate of initial triangular pulse. Left panel Riemann-Liouville model, right panel Caputo model.}
\label{fig2}
\end{figure}

As expected,  the predictions in the left panel of this figure match the results in Fig. 2 in \cite{Baeumer} and the predictions in the right panel match the results in Fig. 5 in \cite{Baeumer} . This verifies that our numerical method, based on the conserved diffusion transport eq. (\ref{con}) and in particular our flux constructions in eq. (\ref{eqrlwflux}) and eq. (\ref{eqcwflux}) are valid.   
 
In terms of the predictions obtained, we note that, at early times the pulse in both models moves to the left. When the pulse encounters the left wall, in both cases the amount of the conserved quantity at the wall increases. Due to the presence of the ``apparent advection'' term, however, we see that, with the RL model the conserved quantity begins to accumulate at the left boundary. At long times both solutions arrive at a steady case. For the C model, in the right panel of Fig. \ref{fig2}, we see that this equilibrium solution is a flat line of unit height. We believe that this is a physically reasonable result because, given sufficient time we should expect a diffusion process, operating in a closed domain, to diffuse any perturbation in the concentration profile.  To the contrary, while the RL  equilibrium result (left panel)  still retains overall conservation, the conserved quantity  clearly accumulates at the left boundary. Effectively, as steady state is approached, the diffusion of the conserved quantity away from the boundary (i.e., to the right) is balanced by the apparent advection of this quantity into the boundary. We do not deny that such a non-local behavior might be physical it is just that our intuition suggests that a ``diffusion'' process should, like the Fourier and C flux laws, always work toward ``flattening''  a potential profile.        
      
Finally, we comment on the zero flux boundary  condition associated with eq. (\ref{eqbal}). If we examine \cite[eq. (8)]{BJDE}, we see that the authors study
the operator $\frac{\partial}{\partial x} \partial^\alpha_{x_C}$, with domain contained in $H^{1+\alpha}$ and they show that it generates a strongly continuous semigroup. Hence, according to Proposition \ref{flux} the zero flux condition is automatically satisfied at $x=0$, if $\alpha>1/2$. However, a corresponding statement for the RL flux is missing in \cite{BJDE}.      

\subsection{Effect of apparent convection with fixed value  boundary conditions}  
Next we consider the following physical problem. Let us assume that $[0,1]$ contains ice at the freezing temperature, $u(x,0)=T_m$ and that this solid phase it maintained in the domain by applying fixed (absorbing) boundary conditions of $u(0,t)=u(1,t)=T_m$. If we were to solve this problem in Warsaw we would probably set $T_m=0^o$C. In this case we can see, by inspection, that any choice of our flux models,  eqs. (\ref{eqgflux}), (\ref{eqrlwflux}), and (\ref{eqcwflux}) in the control-volume discrete form of the conservation balance, (\ref{eqdis}), will, as we might expect, maintain the initial profile through time, $u(x,t)=u(x,0)=0$. For the C flux this is obvious since each local flux contribution to the ``non-local diffusion'' (the sum on the right of eq. (\ref{eqcwflux})) will start and remain at the value of 0. As noted above the RL flux is the C flux plus an advection term, this term vanishes here since the left-hand boundary data is $u(0,t)=T_m=0$. So we see that, when $T_m=0$ the RL flux will also start at and remain at a value of 0. Now, however, let us relocate the solution of the problem to Minneapolis where the freezing temperature value, while just as cold, is given as $T_m=32^o$F. In this location, the Fourier flux and the C flux will still start at and remain at a value of 0. So with these fluxes we will maintain, through time, the solution $u(x,t)=32$, which is a simple translation, according to the temperature conversion factor, of the solution obtained in Warsaw. With the RL flux, however, this will not be the case. Initially, the "non-local diffusion", the first term on the right of (\ref{eqrlwflux}), will take a value of 0 but since it depends directly on the boundary data $u(0,t)=32$, the  ``apparent advection'', the second term on the right of (\ref{eqrlwflux}), will take negative values. The result is that, as time advances, the temperature in the vicinity of the wall will decrease. Computations, indicate that a steady state profile is reached, see Fig. \ref{fig3}, but this is not the expected $u(x,t)=32$, rather a non-horizontal profile forms so that the non-local diffusion contribution to the flux is able to balance the apparent advection contribution. The shape of the resulting profile obtained with the RL flux in Minneapolis cannot be related to the flat profile predicted in Warsaw, this is a physically unrealistic outcome, just {\it rescaling the temperature between Fahrenheit  and Celsius  should not change the shape of the predicted profile}.

\begin{figure}[h]
\centering
\includegraphics[width=6cm]{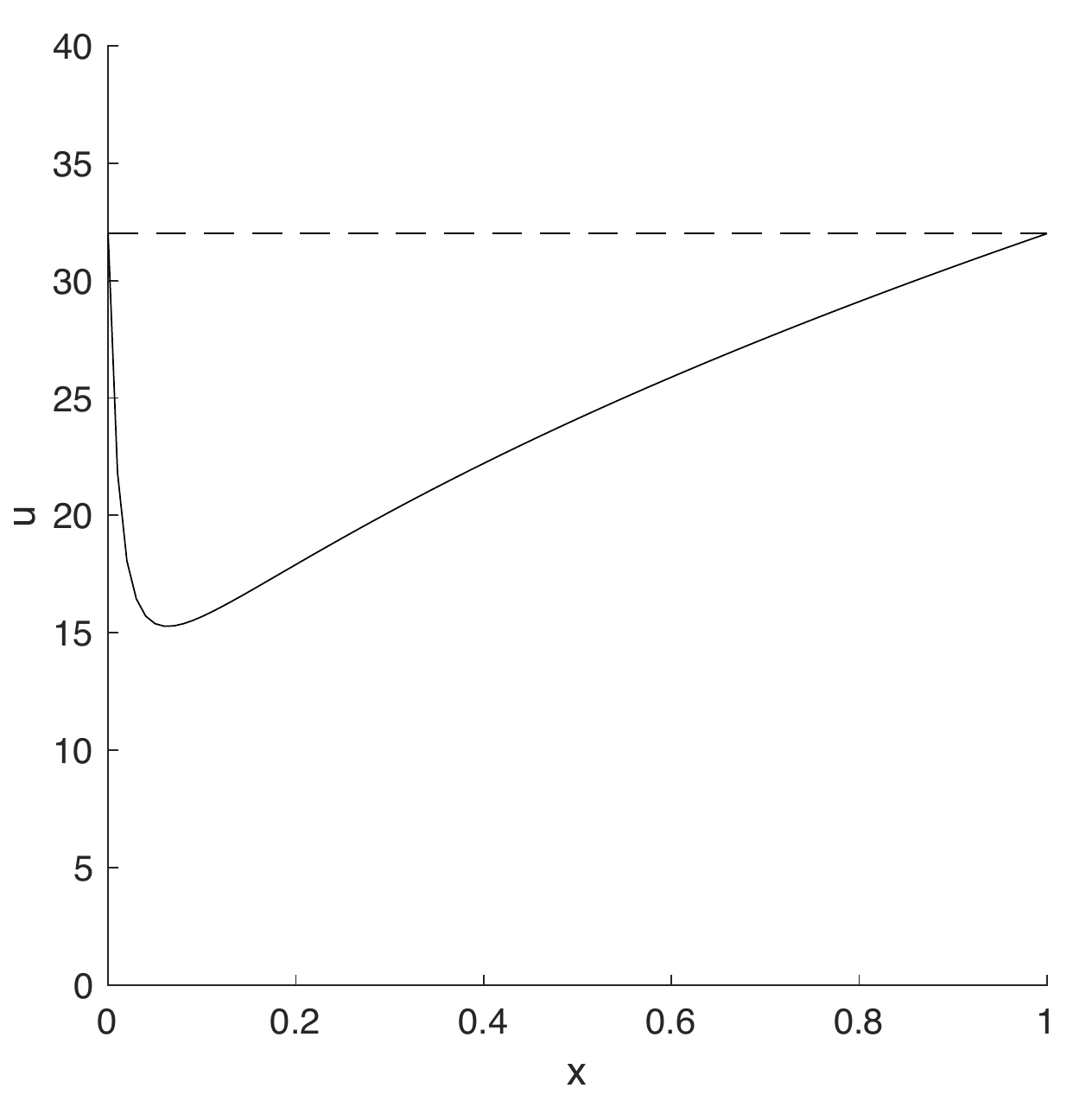}
\caption{Ice domain problem solved ($\Delta x =0.01, \Delta t=0.0005$ ) in  Minneapolis ($T_m=32^o$F) showing close to steady sate solutions using the  C-flux, eq. (\ref{eqcwflux} (dashed line) and RL-flux (eq.(\ref{eqrlwflux} (line).  Note that the presence of the apparent advection term in the RL-flux leads to a non-physical steady sate profile.}
\label{fig3}
\end{figure}

Let us look at this  problem  from the PDE point of view. We consider
\begin{equation}\label{6r1}
 u_t = \partial_{x_{RL}}^{1+\alpha} u\qquad \hbox{in } \Omega\times (0,T),
\end{equation}
because $\partial_{x_{RL}}^{1+\alpha}$ is in the conservative form
or we look at
\begin{equation}\label{6r2}
 u_t = \frac{\partial}{\partial x}\partial_{x_C}^{\alpha} u\qquad \hbox{in } \Omega\times (0,T),
\end{equation}
with non-zero boundary and initial conditions, i.e.
\begin{equation}\label{6r3}
 u|_{\partial \Omega \times(0,T)} = g,\qquad u(x,0) = u_0(x), \quad x\in \Omega.
\end{equation}
Moreover, we assume that the data are consistent,
$$
g|_{t=0} = u_0 |_{\partial\Omega}.
$$
Let us  call by $u^{RL}_g$ (resp. $u^{C}_g$) the solution to (\ref{6r1}) (resp.  
the solution to (\ref{6r2})) with boundary and initial conditions  (\ref{6r3}).

We 
assume that (\ref{6r1}) and (\ref{6r2}) are two variants of the heat  equation. Thermodynamics tells us that the measured temperature does not depend on the scale we use to measure it, either Celsius or Fahrenheit. In other words, if we observe at $t=0$ temperature $u_0$ in Celsius, then the same temperature in Fahrenheit is $\alpha u_0 + \beta$ for certain numbers $\alpha, \beta \in \mathbb{R}$. Of course, we transform the boundary data by the same rule. Moreover, we expect that the temperature at $t>0$ transforms from Celsius to Fahrenheit by the same rule, i.e.
\begin{equation}\label{6r5}
 u^i_{\alpha g + \beta}(t) = \alpha u^i_g (t) + \beta,\qquad i\in\{C, RL\}.
\end{equation}
It turns out that (\ref{6r5}) is valid for $i=C$.

However, (\ref{6r5}) fails for solutions of (\ref{6r1}).
Here is the explanation of this phenomenon. Since (\ref{6r1}) is linear, we have 
$$
u^{RL}_{\alpha g} = \alpha  u^{RL}_g,
$$
but function $v\equiv 1$ is not a steady state solution of (\ref{6r1})
because $\partial_{x_{RL}}^{1+\alpha} 1 = \frac{\alpha(\alpha -1)}{\Gamma(1-\alpha)} x^{-1-\alpha}$. On the
other hand,
$ \frac{\partial}{\partial x}\partial_x^{\alpha} 1 =0$. In other words, changing the temperature scale drastically changes the solution.

\subsection{RL and C flux  calculation of Fig. 7 in \cite{Baeumer}}\label{seven}
We consider Fig. 7 in \cite{Baeumer}, the initial data for problem (\ref{6r1}) or (\ref{6r2}) is
\begin{equation}\label{ccc} 
u(x,0)= \left\{
\begin{aligned}
&\frac{64\pi^3}{\pi^2-4}(x-\tfrac{1}{4})^2\sin(4\pi x) \qquad &\mbox{for }0<x<0.25 \\
&0\qquad &\mbox{ otherwise},
\end{aligned}
\right.
\end{equation}  
and the boundary data is $u(0,t)=u(1,t)=0$. The numerical solution, at times $t=$0.01, 0.04, and 0.2, using the RL flux is shown as the continuous black line in Fig. \ref{fig7}. This solution coincides with the solution using the C flux, red dots. 
\begin{figure}[h]
\centering
\includegraphics[width=6cm]{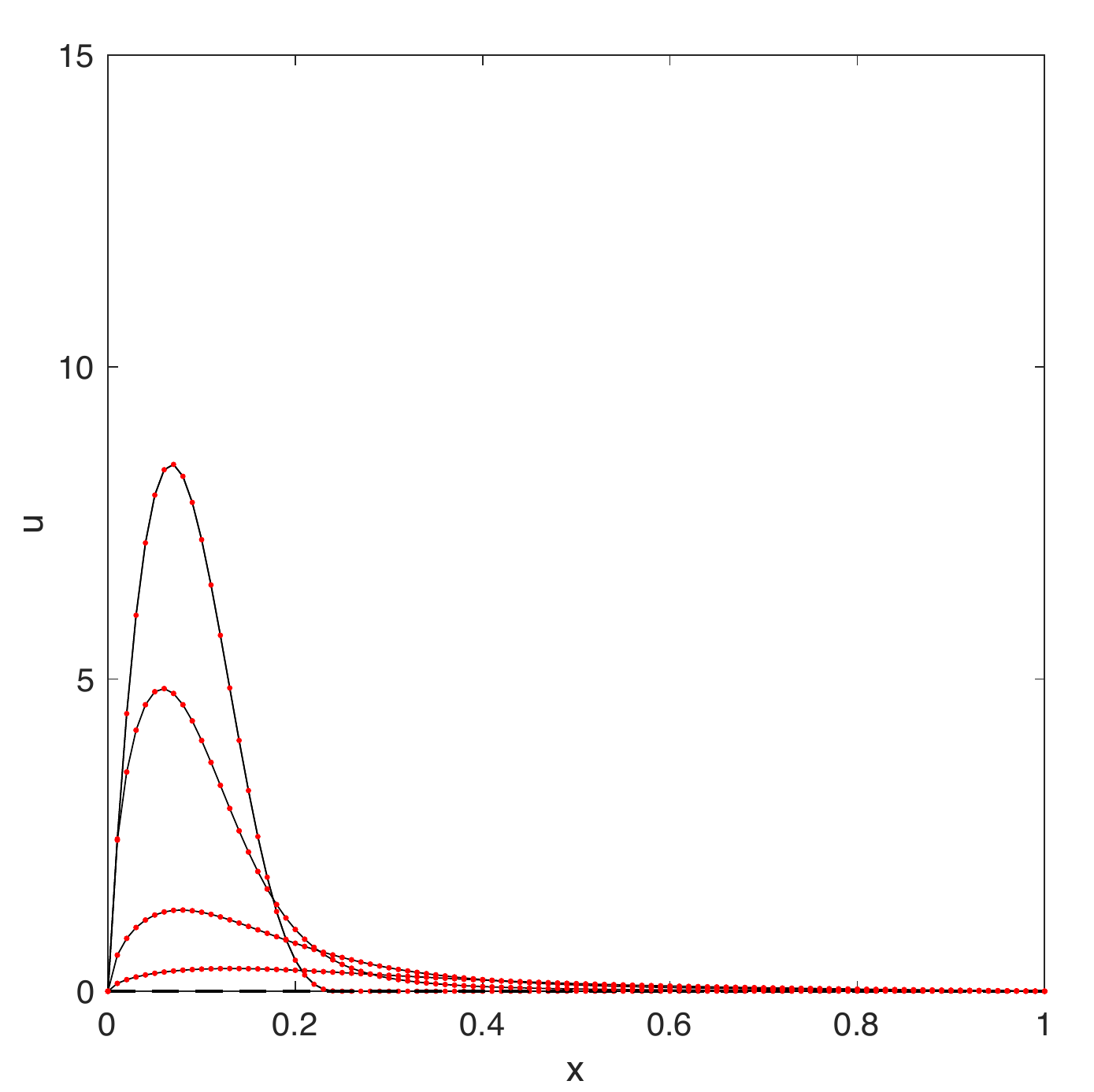}
\caption{Reproduction of Fig 7 with boundary values $u^*(0,t)=u^*(1,t)=0$, the continuous black line is the RL solution, the red points the C solution.}  
\label{fig7}
\end{figure}
\subsection{Fixed non-zero boundary conditions}
We next repeat the calculation but with the modified initial data $u^*=u+5$, and boundary values $u(0,t)=u(1,t)=5$. The C flux solution (red dots), see Fig. \ref{fig7-5} retains the same shape, i.e., the $u$ solution is simply displaced upwards by a constant factor of 5. On the other hand, when we use the RL flux, with this displaced initial data, (continuous line in  Fig. \ref{fig7-5}) we see that we generate transient time step predictions ($t=$0.04 and 0.2) that fall below the initial value (dash line in figure). 
\begin{figure}[h]
\centering
\includegraphics[width=6cm]{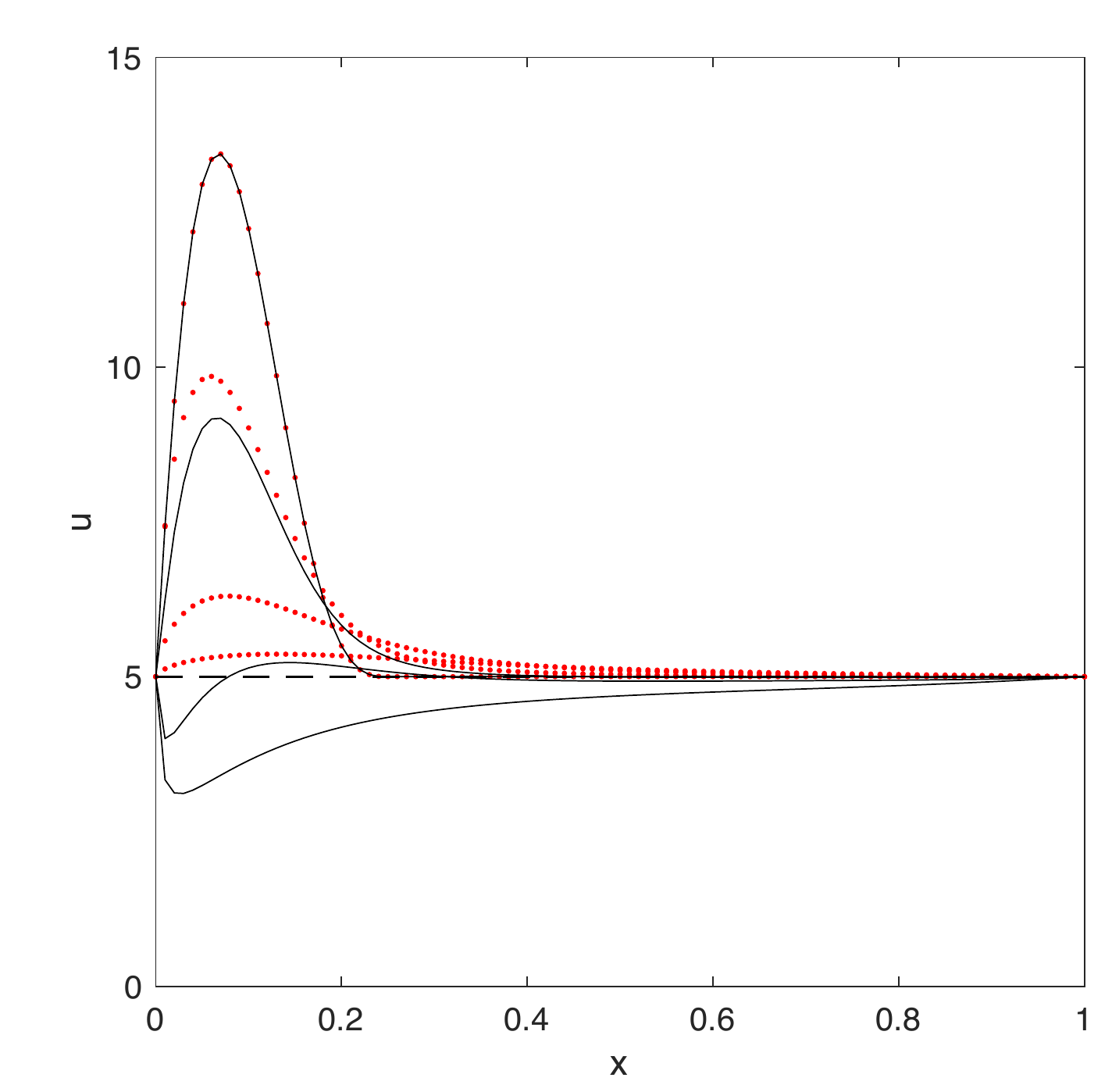}
\caption{Reproduction of Fig 7 with boundary values $u^*(0,t)=u^*(1,t)=5$, the continuous black line is the RL solution, the red points the C solution.}
\label{fig7-5}
\end{figure}
%

This suggest that, with the RL flux, the maximum principle may not hold and physically meaningless predictions might result. Perhaps, meaningful solutions using  the RL flux may only be possible with zero boundary data on the left, i.e., $u(0,t)=0$?

\section{A suggestion for a parsimonious  fractional flux definition} 

From the above it is clear that the differences between conserved non-local flux models can be quite pronounced. This feature significantly undermines our ability to conduct rigorous analysis of the governing equations and to construct physical and consistent numerical solution treatments. One way in which the ambiguity between various non-local flux models might be removed, is to modify the current phenomenological model for the heat flux. Our proposal here is, in an finite interval $x \in [a,b]$, to define the flux at a point $x$ as
\begin{equation}\label{conform}
q(x,t)=\partial_x^{\alpha} \left[ u(x,t)-u(a,t) \right]
\end{equation}   
That is we define the flux as a left-hand fractional space derivative
of $u^*(x,t)=u(x,t)-u(a,t)$. The advantage here is that we do not have
to specifically define the nature of the fractional derivative--since
the C and RL derivatives are equivalent when the initial data is
zero. i.e., when $u^*(x=a, t)=0$, see the Remark at the end of Section
2. Further, at any control volume face in a discretization, the flux can be calculated solely in term of the weighted sum of integer Fourier fluxes that define the non-local diffusion on the right side of eq. (\ref{eqcwflux}) , i.e., any apparent advection contribution to the flux vanishes since, by construction the boundary data $u^*(0,t)=0$.

 \section{Conclusion}
Our work here has focused on the analytical and numerical study of  conserved statements of diffusion transport which employ a non-local  flux-laws  constructed in terms of a fractional derivative operators--space fractional diffusion equations.

The key { contributions include:
   
\begin{itemize}

\item{A presentation and discussion of available results on well posedness of fractional diffusion equations}.

\item{A conserved, control volume discrete solution of fractional diffusion equations, with the development of associated boundary condition treatments that are independent of the specific flux law used.}

\item{The identification of an ``apparent advection'' term in the discrete form of the Riemann-Liouville fractional flux model that can lead to predictions that are may not be physically sound.}

\item Stressing the role of the maximum principle as an important tool in the investigation of solutions. In particular, problems violating this principle may lead to unsound physical models.

\item{The proposal of a  parsimonious fractional flux definition that mitigates the ambiguity associated with the different predictive outcomes of the Riemann-Liouville and Caputo fractional flux definitions.}

\end{itemize}

\section*{Acknowledgment}
The authors thank BIRS in Banff for creating a stimulating environment at the workshop {\it Advanced Developments for Surface and Interface Dynamics - Analysis and Computation}, where this paper was initiated.
VV notes the support of the James L. Record Chair, from
the Department of Civil, Environmental and Geo- Engineering,
University of Minnesota. 
A part of the research for this paper was done at the University of Warsaw, where PR and VV enjoyed a partial support of the National Science Center, Poland, through the grant number
2017/26/M/ST1/00700.
All authors acknowledge stimulating
discussion with Adam Kubica and Katarzyna Ryszewska.


\begin{thebibliography}{99}
 %
\bibitem{adams} R.A.Adams, J.J.F.Fournier,   Sobolev spaces. Second edition.  Elsevier/Academic Press, Amsterdam, 2003
%
\bibitem{Baeumer}
B.Baeumer, M.Kov\'acs, M.M.Meerschaert, H.Sankaranarayanan,
 Boundary conditions for fractional diffusion, {\it J. Comput. Appl. Math.}, {\bf 336}, (2018), 408-424.
 %
\bibitem{Barxive}
B.Baeumer, T.Luks, M.M. Meerschaert, Space-time fractional Dirichlet problems, {\it arXiv:1604.06421}.

\bibitem{BJDE}
B.Baeumer, M.Kov\'acs, M.M.Meerschaert, H.Sankaranarayanan, Fractional partial differential equations with boundary conditions, {\it J. Differential Equations}, {\bf  264} (2018), no. 2, 1377--1410.
%
\bibitem{Barles} 
G.Barles, E.Chasseigne, A.Ciomaga, C.Imbert, 
Lipschitz regularity of solutions for mixed integro-differential equations. 
{\it J. Differential Equations}, {\bf 252} (2012), no. 11, 6012--6060. 
%
\bibitem{Caffarelli}
L.Caffarelli, L.Silvestre, 
Regularity theory for fully nonlinear integro-differential equations, 
{\it Comm. Pure Appl. Math.}, {\bf  62} (2009), no. 5, 597--638.
%
\bibitem{Joe}
J. Fourier, Analytical Theory of Heat, G.E Stechert \& Co, New York, 1878 (translated
by A. Freeman).
%
\bibitem{GigaNamba}
Y.Giga, T.Namba, 
Well-posedness of Hamilton-Jacobi equations with Caputo's time fractional derivative. 
{\it  Comm. Partial Differential Equations}, {\bf  42} (2017), no. 7, 1088--1120. 
%
\bibitem{Yamamoto}
R. Gorenflo, Y. Luchko, M. Yamamoto, Time-fractional diffusion equation in the fractional Sobolev spaces, {\it  Fract. Calc. Appl. Anal.}, {\bf 18} (2015), no. 3, 799--820.
%
\bibitem{Mark}
M.M. Meerschaert, A. Sikorskii, Stochastic Models for Fractional Calculus, vol. 43, De
Gruyter, Berlin, 2011.
%
\bibitem{Metzler}
R. Metzler, J. Klafter, The random Walk's guide to anomalous diffusion: a fractional
dynamics approach, {\it Phys. Rev.}, {\bf 339} (2000), 1-77.
%
\bibitem{Pat}
S.V.Patankar, Numerical Heat Transfer and Fluid Flow, Hemisphere Publishing, 1980.
%
\bibitem{Pod}
I. Podlubny, Fractional Differential Equations: An Introduction to Fractional Derivatives,
Fractional Differential Equations, to Methods of Their Solution and Some of
Their Applications, Academic Press, San Diego, USA, 1999.
%
\bibitem{NaRy}
T.Namba, P.Rybka, On viscosity solutions of space-fractional diffusion equations of Caputo type, {\it in preparation}.
%
\bibitem{Rys} K.Ryszewska,  An analytic semigroup generated by a fractional differential operator, {\it in preparation}.
%
\bibitem{Rys2} K.Ryszewska,  A space - fractional Stefan problem, {\it in preparation}.
%
\bibitem{Rina}
R. Schumer, M.M. Meerschaert, B. Baeumer, Fractional advection dispersion equations
for modeling transport at the earth surface, {\it J. Geophys. Res.}, {\bf 114} (2009) F00A07.

\bibitem{VollerAHT}
V.R.Voller, Anomalous heat transfer: Examples, fundamentals, and fractional calculus models, {\it Advances in Heat Transfer}, {\bf 50}, (2018), 333-378.
%
\bibitem{VollerBook}
V.R. Voller, Basic Control Volume Finite Element Methods for Fluids and Solids, World Scientific, 2009.
%
\bibitem{Zhang}
Y. Zhang, D.A. Benson, M.M. Meerschaert, E.M. La Bolle, Space-fractional
advectione dispersion equations with variable parameters: Diverse formulas, numerical solutions and application to the macro dispersion experiment site data, {\it Water Resources Res.}, {\bf 43}, (2007), W05439.

\bibitem{Dan}
D.P. Zielinski, V.R. Voller, A control volume finite element method with spines for
solutions of fractional heat conduction equations, {\it Numer. Heat Transfer}, {\bf B 70}, (2016), 503-516.


\end{thebibliography}
\end{document}